\documentclass{birkjour}
\usepackage{amsfonts,amssymb,amsmath,amsgen,amsopn,amsbsy}
 \usepackage{amssymb,stmaryrd}
\usepackage{amsmath}
\usepackage{graphicx}
\usepackage{epstopdf}
\usepackage{xcolor}

 \numberwithin{equation}{section}
\newtheorem{theorem}{Theorem}[section]

\newtheorem{lemma}[theorem]{Lemma}

\newtheorem{corollary}[theorem]{Corollary}

\theoremstyle{definition}
\newtheorem{definition}[theorem]{Definition}

\theoremstyle{remark}

\begin{document}

%
%
%
%
%
%
%
%
%

\setlength{\headheight}{26pt}
\title[Closures of Harmonic Bergman–Besov Spaces ]
 {Closures of Harmonic Bergman–Besov Spaces in the Weighted Harmonic Bloch
Spaces on the Unit Ball}

\author[\"{O}mer Faruk Do\u{g}an]{\"{O}mer Faruk Do\u{g}an }

\subjclass{31B05, 46E15, 42B35,26A33, 46E22}

\keywords{Weighted harmonic Bloch Space, Harmonic Bergman-Besov space, Closure}

\date{January 1, 2004}
\address{Department of Mathematics, Tek$\dot{\hbox{\i}}$rda\u{g} \\ Namik Kemal University,
59030 \\ Tek$\dot{\hbox{\i}}$rda\u{g}, Turkey,}

\email{ofdogan@nku.edu.tr}

\begin{abstract}
In this paper, via invertible radial differential operators, we characterize the closures of the harmonic Bergman-Besov Spaces in the weighted harmonic Bloch spaces on the unit ball of $\mathbb{R}^n$ in terms of natural level sets. To this end, we first
show that the harmonic Bergman–Besov space is contained in the weighted harmonic little Bloch space.

\end{abstract}

\maketitle

\section{Introduction}\label{s-introduction}
Let $\mathbb{B}=\mathbb{B}_{n}$ be the open unit ball and $\mathbb{S}=\mathbb{S}_{n}$ be the unit sphere in $\mathbb{R}^n$  for $n\geq 2$.
Let $\nu$  be the normalized Lebesgue measure on $\mathbb{B}$.
For  $\alpha\in \mathbb{R}$, we define the weighted measures $\nu_\alpha$ on $\mathbb{B}$ by
\[
d\nu_\alpha(x)=\frac{1}{V_\alpha} (1-|x|^2)^\alpha d\nu(x).
\]
These measures are finite only when $\alpha>-1$ and in this case we select the constant $V_\alpha$ so that $\nu_\alpha(\mathbb{B})=1$. Naturally $V_0=1$. If $\alpha\leq -1$, we set $V_\alpha=1$.  For $0<p<\infty$, we denote the Lebesgue classes with respect to $\nu_\alpha$ by $L^p_{\alpha}$ and the corresponding norms by $\|\cdot\|_{L^p_{\alpha}}$.

For a multi-index $m=(m_1,\dots,m_n)$ where $m_1,\dots,m_n$ are non-negative integers and for smooth $f$ we write
\[
\partial^m f= \frac{\partial^{|m|} f}{\partial x_1^{m_1}\cdots\partial x_n^{m_n}},
\]
where $|m|=m_1+\dots+m_n$.

Let $h(\mathbb{B})$ be the space of all complex-valued harmonic functions on $\mathbb{B}$ endowed with the topology of uniform convergence  on compact subsets. The space of bounded harmonic functions on $\mathbb{B}$ is denoted by $h^{\infty}$. For $\alpha>-1$, the harmonic Bergman space $b^p_\alpha$ is defined as $b^p_\alpha=  L^p_\alpha \cap h(\mathbb{B})$. This family of spaces can be extended to all $\alpha\in \mathbb{R}$.  For $\alpha\in\mathbb{R}$ and $0<p<\infty$, let $N$ be a non-negative integer such that
\begin{equation}\label{alpha-pN}
\alpha+pN>-1.
\end{equation}
The harmonic Bergman-Besov space $b^p_\alpha$ consists of all $f\in h(\mathbb{B})$ such that
\[
(1-|x|^2)^N \partial^m f \in L^p_\alpha,
\]
for every multi-index $m$ with $|m|=N$.

The space $b^p_\alpha$ do not depend on the choice $N$ as long as (\ref{alpha-pN}) is satisfied. When $\alpha>-1$, one can choose $N=0$ and the resulting space is harmonic weighted Bergman space. When $\alpha=-1$ and $p=2$, the space $b^2_{-1}$ is the harmonic Hardy space. When $\alpha=-n$, the measure $d\nu_{-n}$ is M\"obius invariant and the spaces $b^p_{-n}$ are called by many authors harmonic Besov spaces. If, in addition, $p=2$, the space $b^2_{-n}$ is the harmonic Dirichlet space.

  When $p=2$, the spaces $b^2_\alpha$ are   reproducing kernel Hilbert spaces whose kernels play a major
role in the theory of harmonic  Bergman–Besov spaces $b^p_\alpha$. The spaces $b^p_\alpha$  can  be defined by using the radial differential operators $D^t_s$ $(s,t \in \mathbb{R})$ introduced in \cite{GKU1} and \cite{GKU2}. These operators are defined in terms of reproducing kernels of harmonic Besov spaces and are specific to these spaces, but still mapping  $h(\mathbb{B})$ onto itself.  The properties of $D^t_s$ will be reviewed in Section \ref{s-preliminaries}. Also consider the linear transformation $I_{s}^{t}$ defined for $f\in h(\mathbb{B})$ by
\begin{equation*}
  I^t_s f(x) := (1-|x|^2)^t D^t_s f(x).
\end{equation*}
For $0<p<\infty$ and $\alpha \in \mathbb{R}$, we define the harmonic Bergman-Besov space $b^p_\alpha$ to consist of all $f\in h(\mathbb{B})$ for which $ I^t_s f$
belongs to  $L^p_\alpha$ for some  $s,t$ satisfying (see \cite{GKU2} when $1\leq p<\infty$, and \cite{DOG} when $0<p<1$)
\begin{equation}\label{alpha+pt}
 \alpha+pt>-1.
 \end{equation}
The quantity
\[
\|f\|^p_{b^p_\alpha} = \| I^t_s f\|^p_{L^p_\alpha}=\frac{1}{V_\alpha}\int_{\mathbb{B}} |D^t_s f(x)|^p (1-|x|^2)^{\alpha+pt} d\nu(x) <\infty
\]
defines a norm (quasinorm when $0<p<1$) on $b^p_\alpha$ for any such $s,t$.

It is well-known that the above definition is independent of $s,t$ under (\ref{alpha+pt}),  and the norms (quasinorms when $0<p<1$) on a given space are all equivalent. Thus for a given pair $s,t$,  $ I^t_s$ isometrically imbeds $b^p_\alpha$ into $L^p_\alpha$ if and only if (\ref{alpha+pt}) holds. These matters are studied in detail in \cite{GKU2} when $1\leq p<\infty$ and in \cite{DOG} when $0<p<1$ and will be reviewed in Section \ref{s-preliminaries}.

The well-known harmonic Bloch space $b$ is consists of all $f\in h(\mathbb{B})$ such that
\begin{equation*}
\sup_{x\in \mathbb{B}} (1-|x|^2) |\nabla f(x)|
\end{equation*}
is finite. Let $L^{\infty}$ be the Lebesgue class of essentially bounded functions on $\mathbb{B}$. For $\alpha \in \mathbb{R}$ we define
\begin{equation*}
L^{\infty}_\alpha = \{\varphi : (1-|x|^2)^{\alpha} \varphi(x) \in L^{\infty} \},
\end{equation*}
so that $L^{\infty}_0 = L^{\infty}$ and with norm $\|\varphi\|_{L^{\infty}_\alpha} = \|(1-|x|^2)^\alpha \varphi(x)\|_{L^\infty}$.
For $\alpha>0$, the weighted harmonic Bloch space $b^\infty_\alpha$  is $h(\mathbb{B})\cap L^\infty_\alpha$ endowed with the norm $ \|\cdot\|_{L^{\infty}_\alpha}$, and  is clearly imbedded in  $  L^\infty_\alpha$ by the inclusion $i$.
For  $\alpha \in \mathbb{R}$, we define the weighted harmonic Bloch space $b^\infty_\alpha$ to consist of all $f\in h(\mathbb{B})$ for which $ I^t_s f$ belongs to  $L^\infty_\alpha$ for some  $s$ and $t$ satisfying (see \cite{DU1} )
\begin{equation}\label{alpha+t}
 \alpha+t>0.
 \end{equation}
The quantity
\[
\|f\|_{b^\infty_\alpha} = \| I^t_s f\|_{L^\infty_\alpha}=\sup_{x\in \mathbb{B}}\, (1-|x|^2)^{\alpha+t} |D^t_s f(x)|<\infty,
\]
defines a norm on $b^\infty_\alpha$ for any such $s,t\in \mathbb{R}$.

The weighted harmonic little Bloch space $b^{\infty}_{\alpha 0}$ consist of functions $f\in b^\infty_\alpha$ such that
\[
\lim_{|x|\to 1^{-}} (1-|x|^2)^{\alpha} |I^t_s f(x)|=0,
\]
where $s,t\in \mathbb{R}$ satisfies  $\alpha+t>0$. If $ \alpha=0$, then the space $b_{00}$ is the well-known harmonic little Bloch space.

We mention some properties of the  spaces $b^\infty_{\alpha 0}$ and $b^\infty_{\alpha}$ that can be obtained directly
from their definition. If $f\in h(\overline{\mathbb{B}})$, then $f\in b_{\alpha 0}$; in particular every $b_{\alpha 0}$ (and $b_\alpha$) contains harmonic polynomials and therefore is non-trivial. It is also clear that
\begin{equation*}
  b^{\infty}_\alpha \subset b^{\infty}_{\beta 0} \subset b^{\infty}_\beta \qquad (\text{for\ } \alpha < \beta).
\end{equation*}
The above inclusions are in fact strict (see \cite[Remark 4.9]{DU1} ) and therefore all these spaces are different. It is known that $b_{\alpha 0}$ is closure of polynomials in $b^{\infty}_\alpha$, see \cite[Corollary 4.7]{DU1}.

Note also that above definitions are independent of $s,t$ under (\ref{alpha+t}),  and the norms  in these spaces are all equivalent. Therefore the operator $I^t_s$ isometrically imbeds $b^\infty_\alpha$ into $L^\infty_\alpha$ for a given pair $s,t$ if and only if (\ref{alpha+t}) holds. These matters are studied in detail in \cite{DU1}.

In \cite{DU2},  the precise inclusion relations among  harmonic Bergman-Besov spaces $b^p_\alpha$ and weighted Bloch spaces $b^\infty_\alpha$ on the unit ball of $\mathbb{R}^n$ are given. They showed that if $\alpha \in \mathbb{R}$, then $b^p_{p\alpha -n} \subset b^\infty_\alpha$ and $b^\infty_\alpha$ cannot be replaced by a smaller weighted harmonic Bloch space. But, the inclusion relation between harmonic Bergman–Besov spaces and weighted harmonic little Bloch spaces was not clear from their result, see \cite{DU2}.  The first aim of this paper is to determine the  inclusion relation between these spaces.
\begin{theorem}\label{Theorem-A}
Let  $\alpha \in \mathbb{R}$ and $0<p<\infty$.  Then the harmonic Bergman–Besov space $b^p_{p\alpha-n}$ is  contained in the weighted harmonic little Bloch space $b^\infty_{\alpha 0}$.
\end{theorem}

In $1974$, Anderson et al. \cite{ACP} posed a question that what is the closure of the space of bounded holomorphic functions in the holomorphic Bloch space on the unit disk of complex plane, which is still an open problem.  Later, Ghatage and Zheng \cite{GZ} characterized the
closure of BMOA in the holomorphic Bloch space on the unit disk.  Motivated by these,
there has been numerous subject, namely to describe the closure of several subspaces of the holomorphic Bloch space in the Bloch norm, see \cite{AZ,BG,GMP,LR,MZ,MN,TM,XW,ZR}. Recently, G\"{o}\u{g}\"{u}\c{s} and Yılmaz \cite{GY} characterized the closures of holomorphic Bergman–Besov spaces in weighted holomorphic Bloch spaces.

In this paper, we are going to investigate the closures of harmonic Bergman–Besov spaces in weighted harmonic Bloch spaces.

 If $A$ is a subspace of the weighted harmonic Bloch space $b^\infty_{\alpha}$, then $\mathcal{C}_{b^\infty_{\alpha}}(A)$ will denote the
closure of $A$ in the $b^\infty_{\alpha}$-norm, and the distance from $f\in b^\infty_{\alpha}$ to the subspace $A$ in the
$b^\infty_{\alpha}$-norm will be denoted by  $\text{dist}_{b^\infty_{\alpha}}(f, A)=\inf_{g \in A}\|f-g\|_{b^\infty_{\alpha}}$.

Let $\alpha \in \mathbb{R}$ and $\varepsilon>0$. If $f \in h(\mathbb{B})$, we define the level set $\Omega^{\alpha,t}_{\varepsilon}(f)$ by
\begin{equation*}
\Omega^{\alpha,t}_{\varepsilon}(f) = \{x \in \mathbb{B} : (1-|x|^2)^{\alpha} |I^t_s f(x)|\geq \varepsilon \},
\end{equation*}
where $s, t \in \mathbb{R}$  with $\alpha+t>0$. Let $\chi_{\Omega^{\alpha,t}_{\varepsilon}(f)}$ be the characteristic function of the set $\Omega^{\alpha,t}_{\varepsilon}(f)$.

We can now state our main result.
\begin{theorem}\label{Theorem-1}
Let  $1\leq p<\infty$, $\alpha \in \mathbb{R}$, $\beta\leq p\alpha-n$ and choose $t\in \mathbb{R}$ such that $\alpha+t>0$.
If $f \in b^\infty_{\alpha}$, then the following conditions are equivalent:
\begin{enumerate}
\item[(i)]  $f \in b^\infty_{\alpha 0}$.
\item[(ii)] $f \in \mathcal{C}_{b^\infty_{\alpha}}$ ($b^p_{\beta}$).
\item[(iii)] For every $\varepsilon>0$,
\begin{equation}\label{levelseteq}
\int_{\Omega^{\alpha,t}_{\varepsilon}(f)} (1-|x|^2)^{-n} d\nu(x)< \infty.
\end{equation}
\end{enumerate}
\end{theorem}
If $\beta> p\alpha-n$, then there is no inclusion relation between $b^p_{\beta}$ and $b^\infty_{\alpha}$, hence,
we consider the closure of $b^p_{\beta}\cap b^\infty_{\alpha}$ in the weighted Bloch norm. By \cite[Theorem 1.3]{DU2}, if $\beta>p\alpha-1$, then the closure of $b^p_{\beta}\cap b^\infty_{\alpha}$ in the weighted Bloch norm is $b^\infty_{\alpha}$. It remains to consider the case $ p\alpha-n<\beta \leq p\alpha-1$. In this case, we have the following result.
\begin{theorem}\label{Theorem-2}
Let  $1\leq p<\infty$, $\alpha \in \mathbb{R}$, $ p\alpha-n<\beta \leq p\alpha-1$ and choose $t\in \mathbb{R}$ such that $\alpha+t>0$.
If $f \in b^\infty_{\alpha}$, then the following conditions are equivalent:
\begin{enumerate}
\item[(i)] $f \in \mathcal{C}_{b^\infty_{\alpha}}$ ($b^p_{\beta}\cap b^\infty_{\alpha}$).
\item[(ii)] There exists $t_{0}>t$ with $\alpha+t_{0}>n$ and $\beta+pt_{0}>-1$ such that for every $\varepsilon>0$,
\begin{equation}\label{levelseteq2}
\int_{\Omega^{\alpha,t_{0}}_{\varepsilon}(f)} (1-|x|^2)^{\beta-p\alpha} d\nu(x)< \infty.
\end{equation}
\end{enumerate}
\end{theorem}
The proofs of our results are inspired by the work of G\"{o}\u{g}\"{u}\c{s} and Yılmaz \cite{GY}, where the closures of holomorphic Bergman–Besov spaces in weighted holomorphic Bloch spaces are characterized. 

We briefly summarize some known facts about the harmonic
Bergman-Besov and weighted Bloch spaces in Section  \ref{s-preliminaries}. Section \ref{s-carleson} is devoted to list the some important results that we apply in the proofs. In Section \ref{Inclusion}, we prove Theorem \ref{Theorem-A}. Thus we give the inclusion relation between harmonic Bergman–Besov space and weighted harmonic little Bloch space.  In Section \ref{proof1}, based on the integral representations for the functions in these spaces and some integral estimates involving reproducing kernels, we will give the proof of our main results, Theorems \ref{Theorem-1} and \ref{Theorem-2}.

Throughout the paper, for two positive expressions $A$ and $B$, $A\lesssim B$ means that there exists a positive constant $C$, whose exact value is inessential, such that $A\leq CB$. We also use $A\sim B$ to mean both $A\lesssim B$ and $B\lesssim A$.
\section{PRELIMINARIES}\label{s-preliminaries}
In this section we collect some known facts that will be used later.

The Pochhammer symbol $(a)_b$ is defined by
\[
(a)_b=\frac{\Gamma(a+b)}{\Gamma(a)},
\]
when $a$ and $a+b$ are off the pole set $-\mathbb{N}$ of the gamma function.  Stirling formula provides
\begin{equation}\label{Stirling}
\frac{(a)_c}{(b)_c} \sim c^{a-b} \quad (c\to\infty).
\end{equation}

A harmonic function $f$ on $\mathbb{B}$ has a homogeneous expansion, that is there exist homogeneous harmonic polynomials $f_k$ of degree $k$ such that $f(x)=\sum_{k=0}^\infty f_k(x)$. The series uniformly and absolutely converges on compact subsets of $\mathbb{B}$.

When $\alpha>-1$ the point evaluation functional $f\mapsto f(x)$ is bounded on the Hilbert space $b^2_\alpha$, so by the Riesz representation theorem there exists $R_\alpha(x,y)$ such that
\[
f(x)=\int_{\mathbb{B}} f(y) \overline{R_\alpha(x,y)}\, d\nu_\alpha(y)  \quad (\alpha>-1).
\]
for every $f\in b^2_\alpha$ and $x\in \mathbb{B}$. It is well-known that $R_\alpha$ is real-valued and  symmetric in its variables. The homogeneous expansion of $R_\alpha(x,y)$ can be expressed in terms of zonal harmonics $Z_k(x,y)$ (see \cite{DS}, \cite{M})
\[
R_\alpha(x,y)=\sum_{k=0}^\infty \frac{(1+n/2+\alpha)_k}{(n/2)_k} Z_k(x,y)=:\sum_{k=0}^\infty \gamma_k(\alpha) Z_k(x,y), \qquad (\alpha>-1).
\]
For definition and details about zonal harmonics, see \cite[Chapter 5]{ABR}. We have by (\ref{Stirling})
\begin{equation}\label{gamma_k}
\gamma_k(\alpha)\sim k^{\alpha+1}.
\end{equation}

The reproducing kernels $R_\alpha(x,y)$ can be extended to all $\alpha\in \mathbb{R}$ (see \cite{GKU1,GKU2}), where the crucial point is not the precise form of the kernel but preserving the property (\ref{gamma_k}).

\begin{definition}
Let $\alpha\in\mathbb{R}$. Define
\[
\gamma_k(\alpha):= \begin{cases}
\dfrac{(1+n/2+\alpha)_k}{(n/2)_k}, &\text{if $\, \alpha>-(1+n/2)$}; \\
\noalign{\medskip}
\dfrac{((1)_k)^2}{(1-(n/2+\alpha))_k (n/2)_k}, &\text{if $\, \alpha\leq -(1+n/2)$};
\end{cases}
\]
and $\displaystyle R_\alpha(x,y):=\sum_{k=0}^\infty \gamma_k(\alpha) Z_k(x,y)$.
\end{definition}

For every $\alpha\in \mathbb{R}$,  $R_\alpha(x,y)$ is harmonic as a function of either of its variables on $\overline{\mathbb{B}}$. By (\ref{Stirling}), the property (\ref{gamma_k}) holds for all $\alpha\in\mathbb{R}$. By using the coefficients $ \gamma_k(\alpha)$ in the extended kernels, we define the radial differential operators $D^t_s$ of order $t$.
\begin{definition}
Let $f\in h(\mathbb{B})$ and $f=\sum_{k=0}^\infty f_k$ be its homogeneous expansion. For $s,t\in\mathbb{R}$ we define  $D_s^t$ of order $t$ by
\begin{equation*}
  D_s^t f := \sum_{k=0}^\infty \frac{\gamma_k(s+t)}{\gamma_k(s)} \, f_k.
\end{equation*}
\end{definition}
By (\ref{gamma_k}), $\gamma_k(s+t)/\gamma_k(s) \sim k^t$ for any $s,t$. For every $s\in \mathbb{R}$, $D_s^0=I$, the identity. The additive property $D_{s+t}^{z} D_s^t = D_s^{z+t}$  of $D^t_s$ implies that it is invertible with two-sided inverse
\begin{equation}\label{*}
D^{-t}_{s+t} D^t_s = D^t_s D^{-t}_{s+t} = I.
\end{equation}

For every $s,t \in \mathbb{R}$, the operator $D^t_s: h(\mathbb{B})\to h(\mathbb{B})$ is continuous in the topology
of uniform convergence on compact subsets (see \cite[Theorem 3.2]{GKU2}). The operator $D^t_s$ is constructed so that in all cases
\begin{equation}\label{**}
D_s^t R_s(x,y)=R_{s+t}(x,y),
\end{equation}
where differentiation is performed on one of the variables.

One of the most important properties about the operators $D^t_s$ is that it allows us to pass from one Bergman-Besov (or Bloch) space to another. More precisely, we have the following results.

\begin{lemma}\label{Apply-Dst}
Let $0<p<\infty$ and $\alpha,s,t\in \mathbb{R}$.
\begin{enumerate}
  \item[(i)] The map $D^t_s:b^p_\alpha \to b^p_{\alpha+pt}$ is an isomorphism.
  \item[(ii)] The map $D^t_s:b^\infty_\alpha \to b^\infty_{\alpha+t}$ is an isomorphism.
\end{enumerate}
\end{lemma}

For a proof of part (i) of the above lemma see \cite[Corollary 9.2]{GKU2} when $1\leq p<\infty$ and \cite[Proposition 4.7]{DOG} when $0<p<1$. For part (ii) see \cite[Proposition 4.6]{DU1}.

\section{Main Tools}\label{s-carleson}
In this section, we will give some known lemmas that will be important to the proof our main results.
The next lemma provides an estimate on weighted integrals of powers of $R_\alpha(x,y)$.
When $\alpha>-1$ and $w>0$, it is proved in \cite[Proposition 8]{M}. For the whole range $\alpha\in \mathbb{R}$ see \cite[Theorem 1.5]{GKU2}.

\begin{lemma}\label{norm-kernel}
Let $\alpha\in \mathbb{R}$, $0<p<\infty$ and $d>-1$. Set $w=p(n+\alpha)-(n+d)$. Then
\[
\int_{\mathbb B}|R_\alpha(x,y)|^p\,(1-|y|^2)^d\,d\nu(y)
\sim\begin{cases}
1,&\text{if $w<0$};\\
\noalign{\medskip}
1+\log\dfrac1{1-|x|^2},&\text{if $w=0$};\\
\noalign{\medskip}
\dfrac1{(1-|x|^2)^w},&\text{if $w>0$}.
\end{cases}
\]
\end{lemma}

We require the inclusion relations between harmonic Bergman-Besov and we\-igh\-ted Bloch spaces. We refer to \cite{DU2} for results on inclusions where also references  to earlier work can be found.
We have the following inclusion relations between harmonic Bergman-Besov spaces. For  proofs see \cite[Theorems 1.1 and 1.2]{DU2}.
\begin{theorem}\label{Besov-Besov1}
Let $0<q<p<\infty$ and $\alpha,\beta\in \mathbb{R}$. Then
\[
b^p_\alpha \subset b^q_\beta \quad \text{if and only if} \quad \frac{\alpha+1}{p} < \frac{\beta+1}{q}.
\]
\end{theorem}

\begin{theorem}\label{Besov-Besov2}
Let $0<p\leq q<\infty$ and $\alpha,\beta\in\mathbb{R}$. Then
\[
b^p_\alpha \subset b^q_\beta \quad \text{if and only if} \quad \frac{\alpha+n}{p} \leq \frac{\beta+n}{q}.
\]
\end{theorem}

We also have the following inclusion relation between a Bergman-Besov space $b^p_\alpha$ and a weighted Bloch space $b^\infty_\beta$. For a proof see \cite[Theorem 1.3]{DU2}.
\begin{theorem}\label{Besov-Bloch}
Let $0<p<\infty$ and $\alpha,\beta\in \mathbb{R}$. Then

\begin{enumerate}
  \item[(i)] $b^\infty_\alpha \subset b^p_\beta$ if and only if $\displaystyle \alpha<\frac{\beta+1}{p}$.
  \item[(ii)] $b^p_\beta \subset b^\infty_\alpha$ if and only if $\displaystyle \alpha \geq \frac{\beta+n}{p}$.
\end{enumerate}
\end{theorem}
Note that all the inclusions above are continuous, strict, and the best possible.

Fix $\zeta\in \mathbb{S}$. Then for any $s\in \mathbb{R}$, we have $R_s(\,\cdot\, ,\zeta) \in h(\mathbb{B})$. In the next theorem we determine when $R_s(\, \cdot\, ,\zeta)$ belongs to $b_{\beta}^{p}$ and therefore provide non-trivial (i.e. non-polynomial) examples of elements of $b_{\beta}^{p}$.
\begin{theorem}\label{nonpolex}
Let $1\leq p<\infty$ and $s,\beta \in \mathbb{R}$. Fix $\zeta\in \mathbb{S}$. Then  $R_{s}(\,x\, ,\zeta) \in b_{\beta}^{p}$
if and only if  $\beta+n>p(n+s)$.
\end{theorem}

\begin{proof}
Pick $t\in \mathbb{R}$  large enough that $\beta+pt>-1$ and $n+\beta+t>0$ . By (\ref{**}), we have
\begin{equation*}
I^{t}_{s} R_s(\,x\, ,\zeta)= (1-|x|^{2})^{t}D^{t}_{s} R_s(\,x\, ,\zeta)=(1-|x|^{2})^{t} R_{s+t}(\,x\, ,\zeta).
\end{equation*}
Let $c=p(n+s+t)-(n+\beta+pt)=p(n+s)-(n+\beta)$. If $c<0$, then by Lemma \ref{norm-kernel}, we have
\begin{align*}
\frac{1}{V_{\beta}}\int_{\mathbb{B}} |I_{s}^{t} R_{s}(\,x\, ,\zeta)|^{p} (1-|x|^{2})^{\beta} \, d\nu(x)
&=\frac{1}{V_{\beta}}\int_{\mathbb{B}} | R_{s+t}(\,x\, ,\zeta)|^{p} (1-|x|^{2})^{\beta+pt} \, d\nu(x)\\
&\sim 1.
\end{align*}
 Hence, if $(n+\beta)>p(n+s)$, then  $I_{s}^{t} R_{s}(\,x\, ,\zeta)\in L^{p}_{\beta}$ and therefore $R_{s}(\,x\, ,\zeta)\in b^{p}_{\beta}$.
For the reverse implications, if $(n+\beta)\leq p(n+s)$, then the above integral  diverges by Lemma \ref{norm-kernel}. Thus, for $(n+\beta)\leq p(n+s)$, we have $I_{s}^{t} R_{s}(\,x\, ,\zeta)\notin L^{p}_{\beta}$  and so $R_{s}(\,x\, ,\zeta)\notin b^{p}_{\beta}$.
\end{proof}

We also need the following less known Minkowski integral inequality that in effect exchanges the order of integration; for a proof, see  \cite[Theorem 3.3.5]{O1} for example.

\begin{lemma}\label{Minkowski int-inequality}
Let $(X,A,\mu)$ and $(Y,B,\upsilon)$ be $\sigma$-finite measure
spaces. If $1\leq p\leq \infty$ and   $f(x,y)$ is a measurable function on $X\times Y$, then
\begin{equation*}
\left(\int_{Y}\left(\int_{X} |f(x,y)| d\mu(x)\right)^{p}  d\upsilon(y)\right)^{1/p}\leq \int_{X}\left(\int_{Y} |f(x,y)|^{p} d\upsilon(y)\right)^{1/p}  d\mu(x),
\end{equation*}
with an appropriate interpretation with the $L^{\infty}$ norm when $p=\infty$.
\end{lemma}

\section{Proof of Theorem \ref{Theorem-A}}\label{Inclusion}

In this section, we will prove Theorem \ref{Theorem-A}. 

\begin{proof}[Proof of Theorem \ref{Theorem-A}]
By Theorem \ref{Apply-Dst}, it is enough to show that $b^{p}_{p\alpha-n} \subset b^{\infty}_{\alpha 0}$ for $ \alpha=0$.
Let $f\in b^{p}_{-n}$. Then for any real number $s$, $D^{1}_{s}f  \in b^{p}_{p-n}$. It follows from
\cite[Proposition 13.1]{GKU2} when $0\leq p<\infty$ and \cite[Proposition 5.3]{DOG2}  when $0< p<1$ that
\begin{equation*}
 \lim_{|x|\to 1^{-}} (1-|x|^{2})\big|D^{1}_{s}f (x)\big| =0.
\end{equation*}
Thus $f \in b^{\infty}_{00}$ and $b^{p}_{-n}\subseteq b^{\infty}_{00}$. That finishes the proof.
\end{proof}

\section{Proof of Theorems \ref{Theorem-1} and \ref{Theorem-2}} \label{proof1}
Our goal in this section is to prove Theorems \ref{Theorem-1} and \ref{Theorem-2}. For $\alpha,\beta \in \mathbb{R}$ with $\beta\leq p\alpha-n$, by Theorem \ref{Besov-Bloch}, $b^p_\beta \subset b^\infty_\alpha$. Thus we can consider the closure of $b^p_\beta$ in the weighted harmonic Bloch space $b^\infty_\alpha$. With all these preliminary works, we have laid the groundwork for proving our main results.
\begin{proof}[Proof of Theorem \ref{Theorem-1}]
Theorem \ref{Besov-Besov2} implies that $b^p_\beta \subset b^p_{p\alpha-n}$ whenever $\beta\leq p\alpha-n$.
Hence it is enough to prove theorem only for  $\beta:=p\alpha-n$.
Since the weighted harmonic little Bloch space $b^\infty_{\alpha 0}$ is the closure of the set of harmonic polynomials in the weighted
harmonic Bloch space $b^\infty_\alpha$  and $b^p_\beta$  contains all harmonic polynomials, we obtain that the closure
of $b^p_\beta$ in  $b^\infty_\alpha$ contains $b^\infty_{\alpha 0}$. On the other hand, by Theorem \ref{Theorem-A}, $b^p_\beta \subset b^\infty_{\alpha 0}$. It is clear that the closure of $b^p_\beta $ in $b^\infty_\alpha$ is contained in $b^\infty_{\alpha 0}$. Thus $b^\infty_{\alpha 0}$ equals to the closure of $b^p_\beta $ in $b^\infty_\alpha$, and so statements (i) and  (ii) are equivalent.

(ii) Implies (iii). Let  $f \in \mathcal{C}_{b^\infty_{\alpha}} (b^p_{\beta})$. Then for any $\varepsilon>0$, there exists $g\in b^p_{\beta}$ such that $\|f-g\|_{b^\infty_{\alpha}}\leq \frac{\varepsilon}{2}$ for some $s,t' \in \mathbb{R}$ satisfying $\alpha+t'>0$. By assumption $\alpha+t>0$, we can take $t'=t$. Given such $t \in \mathbb{R}$, let $p_{1}>p$ be such that $p_{1}>\frac{n-1}{\alpha+t}$. Set $\beta_{1}=\alpha p_{1}-n$. Then $\beta_{1}+p_{1}t>-1$ and by Theorem \ref{Besov-Besov2}, we have $b^p_{\beta}\subset b^{p_{1}}_{\beta_{1}}$. Note that we have
\begin{equation*}
 (1-|x|^2)^{\alpha} |I^t_s f(x)|\leq (1-|x|^2)^{\alpha} |I^t_s (f(x)-g(x))|+(1-|x|^2)^{\alpha} |I^t_s g(x)|, \quad x\in \mathbb{B}
\end{equation*}
for such $s,t\in \mathbb{R}$ with $\alpha+t>0$. This gives
\begin{equation*}
\Omega^{\alpha,t}_{\varepsilon}(f)\subseteq \Omega^{\alpha,t}_{\frac{\varepsilon}{2}}(g).
\end{equation*}
Therefore, since $g\in b^p_{\beta}\subset b^{p_{1}}_{\beta_{1}}$, we have
\begin{align*}
\infty &> \int_{\mathbb{B}} |I^t_s g(x)|^{p_{1}}(1-|x|^2)^{\beta_{1}} d\nu(x) \\
&\geq  \int_{\Omega^{\alpha,t}_{\frac{\varepsilon}{2}}(g)} [(1-|x|^2)^{\alpha} |I^t_s g(x)|]^{p_{1}}(1-|x|^2)^{\beta_{1}-\alpha p_{1}} d\nu(x) \\
&\geq \left(\frac{\varepsilon}{2}\right)^{p_{1}} \int_{\Omega^{\alpha,t}_{\frac{\varepsilon}{2}}(g)} (1-|x|^2)^{\beta_{1}-\alpha p_{1}} d\nu(x) \\
&\geq \left(\frac{\varepsilon}{2}\right)^{p_{1}} \int_{\Omega^{\alpha,t}_{\varepsilon}(f)} (1-|x|^2)^{\beta_{1}-\alpha p_{1}} d\nu(x).
\end{align*}
Hence,
\begin{equation*}
\int_{\Omega^{\alpha,t}_{\varepsilon}(f)} (1-|x|^2)^{-n} d\nu(x)< \infty.
\end{equation*}

(iii) Implies (ii). Conversely, fix $\varepsilon>0$ and assume that $f \in b^\infty_{\alpha}$ satisfy (\ref{levelseteq}).  Since $\alpha+t>0$, choose $s\in \mathbb{R}$ such that $s>\alpha-1$. Then by $(7)$ in \cite{DU1} for $f\in b_{\alpha}$, we have the following integral representation
\begin{equation}\label{Reproducing}
f(x)=\frac{1}{V_{s+t}}\int_{\mathbb{B}} R_s(x,y) (1-|y|^2)^{s+t} D^t_s f(y) \, d\nu(y).
\end{equation}
We split the above integral into two parts. Let $f (z)= f_{1}(z)+f_{2}(z)$, where
\begin{equation*}
f_{1}(x)=\frac{1}{V_{s+t}}\int_{\Omega^{\alpha,t}_{\varepsilon}(f)} R_s(x,y) (1-|y|^2)^{s+t} D^t_s f(y) \, d\nu(y).
\end{equation*}
and
\begin{equation*}
f_{2}(x)=\frac{1}{V_{s+t}}\int_{\mathbb{B}\backslash\Omega^{\alpha,t}_{\varepsilon}(f)} R_s(x,y) (1-|y|^2)^{s+t} D^t_s f(y) \, d\nu(y).
\end{equation*}
Under the condition $\alpha+t> 0$, the norms on $b^\infty_{\alpha}$ are all equivalent. So we can take $s, t$ for which
(\ref{Reproducing}) holds. Then, noting that $D_s^t R_s(x,y)=R_{s+t}(x,y)$ and differentiating the above formula of $f_{2}$
under the integral sign (see  \cite[Corollary .2.5]{DU1}), we obtain that
\begin{align*}
|D_s^t f_{2}(x)|&=\left|\frac{1}{V_{s+t}}\int_{\mathbb{B}\backslash\Omega^{\alpha,t}_{\varepsilon}(f)} D_s^t R_s(x,y) (1-|y|^2)^{s+t} D^t_s f(y) \, d\nu(y)\right| \\
&=\left|\frac{1}{V_{s+t}}\int_{\mathbb{B}\backslash\Omega^{\alpha,t}_{\varepsilon}(f)}  R_{s+t}(x,y) (1-|y|^2)^{s+t} D^t_s f(y) \, d\nu(y)\right| \\
&\leq\varepsilon\frac{1}{V_{s+t}}\int_{\mathbb{B}}  |R_{s+t}(x,y)| (1-|y|^2)^{s-\alpha}  \, d\nu(y).
\end{align*}
 Since $s-\alpha>-1$ and $\alpha+t>0$, it follows from Lemma \ref{norm-kernel} that
 \begin{equation*}
|D_s^t f_{2}(x)|\lesssim \varepsilon\frac{1}{V_{s+t}}\frac{1}{(1-|y|^2)^{\alpha+t} }.
\end{equation*}
 This implies that $f_{2} \in b^\infty_{\alpha}$ and $\|f_{2}\|_{b^\infty_{\alpha}}\lesssim \varepsilon\frac{1}{V_{s+t}} $. Consequently, $f_{1} \in b^\infty_{\alpha}$ and $$\|f-f_{1}\|_{b^\infty_{\alpha}}\lesssim \varepsilon\frac{1}{V_{s+t}}.$$
Since $\varepsilon>0$ is arbitrary, the desired result will be established if we show that $f_{1} \in b^p_{\beta}$. Let $t'$ and $s$ satisfying $\beta+pt'>-1$ and $s>\alpha-1$ be given. Differentiating the above formula of $f_{1}$
under the integral sign (see the similar argument in \cite[Corollary .2.5]{DU1}), we obtain
\begin{align*}
D_s^{t'} f_{1}(x)&=\frac{1}{V_{s+t}}\int_{\Omega^{\alpha,t}_{\varepsilon}(f)} D_s^{t'} R_s(x,y) (1-|y|^2)^{s+t} D^t_s f(y) \, d\nu(y) \\
&=\frac{1}{V_{s+t}}\int_{\Omega^{\alpha,t}_{\varepsilon}(f)}  R_{s+t'}(x,y) (1-|y|^2)^{s+t} D^t_s f(y) \, d\nu(y).
\end{align*}
Further, by the definition of $b^{\infty}_{\alpha}$ and Lemma \ref{Minkowski int-inequality},
\begin{align*}
 \int_{\mathbb{B}}&|D_s^{t'} f_{1}(x)|^{p}(1-|x|^2)^{\beta+pt'}\, d\nu(x)
 \leq\frac{1}{(V_{s+t})^{p}}\int_{\mathbb{B}} \bigg[\int_{\Omega^{\alpha,t}_{\varepsilon}(f)} \\
 &|R_{s+t'}(x,y)| (1-|y|^2)^{s+t} |D^t_s f(y)| \, d\nu(y)\bigg]^{p} (1-|z|^2)^{\beta+pt'} \, d\nu(x)\\
 &\leq\frac{1}{(V_{s+t})^{p}}\bigg[\int_{\Omega^{\alpha,t}_{\varepsilon}(f)}
 \bigg( \int_{\mathbb{B}} |R_{s+t'}(x,y)|^{p} (1-|y|^2)^{p(s+t)}\\
 & |D^t_s f(y)|^{p}(1-|x|^2)^{\beta+pt'} \, d\nu(x)\bigg)^{1/p} \, d\nu(y)\bigg]^{p}\\
 &\leq\frac{1}{(V_{s+t})^{p}}\|f\|^{p}_{b^{\infty}_{\alpha}}\bigg[\int_{\Omega^{\alpha,t}_{\varepsilon}(f)}\bigg( \int_{\mathbb{B}}|R_{s+t'}(x,y)|^{p} (1-|y|^2)^{p(s-\alpha)}\\
 &(1-|x|^2)^{\beta+pt'} \, d\nu(x)\bigg)^{1/p} \, d\nu(y)\bigg]^{p}\\
 &=\frac{1}{(V_{s+t})^{p}}\|f\|^{p}_{b^{\infty}_{\alpha}}\bigg[\int_{\Omega^{\alpha,t}_{\varepsilon}(f)}  (1-|y|^2)^{s-\alpha} \\
 &\bigg(\int_{\mathbb{B}}  |R_{s+t'}(x,y)|^{p}  (1-|x|^2)^{\beta+pt'} \, d\nu(x)\bigg)^{1/p} \, d\nu(y)\bigg]^{p}
\end{align*}
In order to apply Lemma \ref{norm-kernel} to the integral
 \begin{equation*}
\int_{\mathbb{B}}  |R_{s+t'}(x,y)|^{p}  (1-|x|^2)^{\beta+pt'} \, d\nu(x)
\end{equation*}
We write $c=p(n+s+t')-(n+\beta+pt')=p(n+s)-(n+\beta)$. Now recall that $\beta=p\alpha-n$. Hence $c=p(n+s-\alpha)$. Since $\beta+pt'>-1$ and $s-\alpha>-1$, we can apply Lemma \ref{norm-kernel} with $c=p(n+s-\alpha)>0$. Thus
 \begin{equation*}
\int_{\mathbb{B}}  |R_{s+t'}(x,y)|^{p}  (1-|x|^2)^{\beta+pt'} \, d\nu(x)\lesssim \frac{1}{(1-|y|^2)^{p(n+s-\alpha)}}.
\end{equation*}
Therefore,
\begin{align*}
 \int_{\mathbb{B}}&|D_s^{t'} f_{1}(x)|^{p}(1-|x|^2)^{\beta+pt'}\, d\nu(z)\\
 &\lesssim\frac{1}{(V_{s+t})^{p}}\|f\|^{p}_{b^{\infty}_{\alpha}}\bigg[\int_{\Omega^{\alpha,t}_{\varepsilon}(f)}  \frac{(1-|y|^2)^{s-\alpha} }{(1-|y|^2)^{n+s-\alpha}} \, d\nu(y)\bigg]^{p}\\
 &=\frac{1}{(V_{s+t})^{p}}\|f\|^{p}_{b^{\infty}_{\alpha}}\bigg[\int_{\Omega^{\alpha,t}_{\varepsilon}(f)}  (1-|y|^2)^{-n} \, d\nu(y)\bigg]^{p}.
\end{align*}
In view of (\ref{levelseteq}), we show that $f_{1} \in b^p_{\beta}$ and the desired result is established.
\end{proof}

Next we prove Theorem \ref{Theorem-2}. The proof of this theorem follows in a manner similar to the proof of Theorem \ref{Theorem-1}. We include a proof for completeness and to illustrate how it follows from the sufficiently higher order derivatives.

\begin{proof}[Proof of Theorem \ref{Theorem-2}]
(i) Implies (ii). Let  $f \in \mathcal{C}_{b^\infty_{\alpha}} (b^p_{\beta}\cap b^\infty_{\alpha})$. Then for any $\varepsilon>0$, there exists $g\in b^p_{\beta}\cap b^\infty_{\alpha}$ such that $\|f-g\|_{b^\infty_{\alpha}}\leq \frac{\varepsilon}{2}$ for some $s,t' \in \mathbb{R}$ satisfying $\alpha+t'>0$. In particular, we can take $t'\geq t+\frac{n}{p}$.  Then for such $t'$, as is mentioned before one can show that
$\Omega^{\alpha,t}_{\varepsilon}(f)\subseteq \Omega^{\alpha,t}_{\frac{\varepsilon}{2}}(g)$. Note also  that, $\beta+pt'>p\alpha-n+pt+n>-1$.
Therefore, since $g\in b^p_{\beta}$, we have
\begin{align*}
\infty &> \int_{\mathbb{B}} |I^t_s g(x)|^{p}(1-|x|^2)^{\beta} d\nu(x) \\
&\geq  \int_{\Omega^{\alpha,t'}_{\frac{\varepsilon}{2}}(g)} [(1-|x|^2)^{\alpha} |I^t_s g(x)|]^{p}(1-|x|^2)^{\beta-\alpha p} d\nu(x) \\
&\geq \left(\frac{\varepsilon}{2}\right)^{p} \int_{\Omega^{\alpha,t'}_{\frac{\varepsilon}{2}}(g)} (1-|x|^2)^{\beta-\alpha p} d\nu(x) \\
&\geq \left(\frac{\varepsilon}{2}\right)^{p} \int_{\Omega^{\alpha,t'}_{\varepsilon}(f)} (1-|x|^2)^{\beta-\alpha p} d\nu(x),
\end{align*}
which establishes (\ref{levelseteq2}).

(iii) Implies (ii). Fix $\varepsilon>0$ and assume that $f \in b^\infty_{\alpha}$ satisfy (\ref{levelseteq2}) for some $t'>t$, where $\alpha+t'>n$ and $\beta+pt'>-1$  hold.  Since $\alpha+t'>0$, choose $s\in \mathbb{R}$ such that $s>\alpha-1$. Define
\begin{equation*}
f_{1}(x)=\frac{1}{V_{s+t'}}\int_{\Omega^{\alpha,t'}_{\varepsilon}(f)} R_s(x,y) (1-|y|^2)^{s+t'} D^{t'}_s f(y) \, d\nu(y).
\end{equation*}
and
\begin{equation*}
f_{2}(x)=\frac{1}{V_{s+t'}}\int_{\mathbb{B}\backslash\Omega^{\alpha,t'}_{\varepsilon}(f)} R_s(x,y) (1-|y|^2)^{s+t'} D^{t'}_s f(y) \, d\nu(y).
\end{equation*}
 By $(7)$ in \cite{DU1} again, for $f\in b_{\alpha}$, we have $f (z)= f_{1}(z)+f_{2}(z)$. As before the proof will be done when we show that $\|f_{2}\|_{b^\infty_{\alpha}}\lesssim \varepsilon$ and $f_{1} \in b^p_{\beta}\cap b^\infty_{\alpha}$. Pick any $s,t$ such that
 $\alpha+t>0$ and $s>\alpha-1$ holds. Then, by differentiating the above formula of $f_{2}$
under the integral sign (see \cite[Corollary 2.5]{DU1}), (\ref{**}) and Lemma \ref{norm-kernel},
\begin{align*}
|D_s^t f_{2}(x)|&=\left|\frac{1}{V_{s+t'}}\int_{\mathbb{B}\backslash\Omega^{\alpha,t'}_{\varepsilon}(f)} D_s^t R_s(x,y) (1-|y|^2)^{s+t'} D^{t'}_s f(y) \, d\nu(y)\right| \\
&=\left|\frac{1}{V_{s+t'}}\int_{\mathbb{B}\backslash\Omega^{\alpha,t'}_{\varepsilon}(f)}  R_{s+t}(x,y) (1-|y|^2)^{s+t'} D^{t'}_s f(y) \, d\nu(y)\right| \\
&\leq\varepsilon\frac{1}{V_{s+t'}}\int_{\mathbb{B}}  |R_{s+t}(x,y)| (1-|y|^2)^{s-\alpha}  \, d\nu(y)\\
&\lesssim \varepsilon\frac{1}{V_{s+t'}}\frac{1}{(1-|y|^2)^{\alpha+t}}.
\end{align*}
 This implies that $f_{2} \in b^\infty_{\alpha}$ and $\|f_{2}\|_{b^\infty_{\alpha}}\lesssim \varepsilon\frac{1}{V_{s+t}} $. Consequently, $f_{1} \in b^\infty_{\alpha}$ and $\|f-f_{1}\|_{b^\infty_{\alpha}}\lesssim \varepsilon\frac{1}{V_{s+t}}.$

We now show that $f_{1} \in b^p_{\beta}$. Using the facts $f,f_{1} \in b^\infty_{\alpha}$, differentiating the above formula of $f_{1}$
under the integral sign and Fubini theorem, we obtain that
\begin{align*}
 \int_{\mathbb{B}}&|D_s^{t'} f_{1}(x)|^{p}(1-|x|^2)^{\beta+pt'}\, d\nu(x)\\
 &\leq \|f_{1}\|^{p-1}_{b^\infty_{\alpha}} \int_{\mathbb{B}}|D_s^{t'} f_{1}(x)|(1-|x|^2)^{\beta-(p-1)\alpha+t'}\, d\nu(x)\\
  &\leq \|f_{1}\|^{p-1}_{b^\infty_{\alpha}} \frac{1}{V_{s+t'}}\int_{\mathbb{B}} (1-|x|^2)^{\beta-(p-1)\alpha+t'}\, d\nu(x)
  \int_{\Omega^{\alpha,t'}_{\varepsilon}(f)} |R_{s+t'}(x,y)| \\
  &(1-|y|^2)^{s+t'}|D^{t'}_s f(y)| \, d\nu(y)\\
&\leq \|f_{1}\|^{p-1}_{b^\infty_{\alpha}}\frac{1}{V_{s+t'}} \int_{\Omega^{\alpha,t'}_{\varepsilon}(f)} (1-|y|^2)^{s+t'}|D^{t'}_s f(y)| \, d\nu(y)\int_{\mathbb{B}} |R_{s+t'}(x,y)|  \\
&(1-|z|^2)^{\beta-(p-1)\alpha+t'}\, d\nu(x).\\
&\leq \|f_{1}\|^{p-1}_{b^\infty_{\alpha}} \|f\|_{b^\infty_{\alpha}}\frac{1}{V_{s+t'}} \int_{\Omega^{\alpha,t'}_{\varepsilon}(f)} (1-|y|^2)^{s-\alpha} \, d\nu(y)\int_{\mathbb{B}} |R_{s+t'}(x,y)|  \\
&(1-|z|^2)^{\beta-(p-1)\alpha+t'}\, d\nu(x).\\
&\leq \|f_{1}\|^{p-1}_{b^\infty_{\alpha}} \|f\|_{b^\infty_{\alpha}}\frac{1}{V_{s+t'}} \int_{\Omega^{\alpha,t'}_{\varepsilon}(f)} (1-|y|^2)^{\beta-p\alpha} \, d\nu(y)\\
\end{align*}
where the last inequality follows from Lemma \ref{norm-kernel} with $d=\beta-(p-1)\alpha+t'>-n+n>-1$ and $ c=n+s+t'-(n+\beta-(p-1)\alpha+t')=s-\alpha+p\alpha-\beta>-1+1>0$.
In view of (\ref{levelseteq2}), we show that $f_{1} \in b^p_{\beta}$ and the desired result is established.
\end{proof}

We now record  a few immediate consequences of  Theorems \ref{Theorem-1} and \ref{Theorem-2}.
\begin{corollary}\label{Corollary-1}
Let  $1\leq p<\infty$, $\alpha \in \mathbb{R}$, $\beta\leq p\alpha-n$ and choose $t\in \mathbb{R}$ such that $\alpha+t>0$.
If $f \in b^\infty_{\alpha}$, then the following conditions are equivalent in the sense of $\sim$:
\begin{enumerate}
\item[(i)]  $\text{dist}_{b^\infty_{\alpha}}(f, b^\infty_{\alpha0})$,
\item[(ii)] $\text{dist}_{b^\infty_{\alpha}}(f, b^p_{\beta})$,
\item[(iii)] $\inf\{\varepsilon:\frac{\chi_{\Omega^{\alpha,t}_{\varepsilon}(f)}(x)}{(1-|x|^2)^{n} }d\nu(x)\text{is a finite measure}\}$.
\end{enumerate}
\end{corollary}
\begin{corollary}\label{Corollary-2}
Let  $1\leq p_{0}<p_{1}<\infty$ and $\alpha \in \mathbb{R}$. Then $$\text{dist}_{b^\infty_{\alpha}}(f, b^{p_{0}}_{p_{0}\alpha-n})=\text{dist}_{b^\infty_{\alpha}}(f, b^{p_{1}}_{p_{1}\alpha-n}).$$
\end{corollary}
\begin{corollary}\label{Corollary-3}
Let  $1\leq p<\infty$ and $\alpha,\beta \in \mathbb{R}$. If $p\alpha-n<\beta\leq p\alpha-1$, then $ b^\infty_{\alpha 0}\subsetneqq
\mathcal{C}_{b^\infty_{\alpha}}$ ($b^p_{\beta}\cap b^\infty_{\alpha}$).
\end{corollary}
\begin{proof}
Since $b^\infty_{\alpha 0}$ is the closure of the set of harmonic polynomials in the weighted
harmonic Bloch space $b^\infty_\alpha$  and $b^p_\beta$  contains all harmonic polynomials, we see that $b^\infty_{\alpha 0}\subseteqq
\mathcal{C}_{b^\infty_{\alpha}}$ ($b^p_{\beta}\cap b^\infty_{\alpha}$).
Fix $\zeta\in \mathbb{S}$. First, for any $\alpha\in \mathbb{R}$, we have $R_{\alpha-n}(\,\cdot\, ,\zeta) \in b^\infty_{\alpha}$ but $R_{\alpha-n}(\,\cdot\, ,\zeta) \notin b^\infty_{\alpha 0}$ by \cite[Remark 4.9]{DU1}. Next, since $\beta>p\alpha-n$, $R_{\alpha-n}(\,\cdot\, ,\zeta) \in b^p_{\beta}$ by Theorem \ref{nonpolex}. Thus  $R_{\alpha-n}(\,\cdot\, ,\zeta) \in b^p_{\beta}\cap b^\infty_{\alpha} \setminus b^\infty_{\alpha 0}$. We complete the proof.
\end{proof}




\begin{thebibliography}{1}
\bibitem{ACP}
Anderson, J.M., Clunie, J., Pommerenke, Ch.: On Bloch functions and normal functions. J. Reine
Angew. Math. 270, 12--37 (1974)
%
\bibitem{AZ}
Aulaskari, R., Zhao, R.: Composition operators and closures of some Möbius invariant spaces in the
Bloch space. Math. Scand. 107(1), 139-–149 (2010)
%
\bibitem{ABR}
Axler, S.,  Bourdon, P.,  Ramey, W., 2001, Harmonic function theory,
2\textit{nd} ed., Grad. Texts in Math., vol. 137, Springer: New York.
%
\bibitem{BG}
Bao, G., G\"{o}\u{g}\"{u}\c{s}, N.G.: On the closures of Dirichlet type spaces in the Bloch space. Complex Anal.
Oper. Theory 13(1), 45--59 (2019)
%
\bibitem{DOG}
Do\u{g}an, \"{O}. F.: Harmonic Besov spaces with small exponents,
Complex Variables and Elliptic Equations 65 (6), 1051--1075 (2020)
%
\bibitem{DOG2}
Do\u{g}an, \"{O}. F.: Positive Toeplitz operators from a harmonic Bergman–Besov space into another,
Banach J. Math. Anal., 16(70), (2022)
%
\bibitem{DU1}
Do\u{g}an,  \"{O}. F.,  \"Ureyen, A. E.: Weighted harmonic Bloch spaces on the ball,
Complex Anal. Oper. Theory, 12 (5), 1143--1177 (2018)
%
\bibitem{DU2}
 Do\u{g}an, \"{O}. F., \"Ureyen A. E.: Inclusion relations between harmonic Bergman-Besov and weighted Bloch spaces on the
unit ball, Czechoslovak Mathematical Journal, 69(2), 503--523 (2019)
%
\bibitem{DS}
Djrbashian, A. E.,   Shamoian, F. A.: 1988, Topics in the theory of $A^p_\alpha$ spaces,
Teubner Texts in Mathematics: 105, BSB B. G. Teubner Verlagsgesellschaft, Leipzig.
%
\bibitem{GMP}
Galanopoulos, P., Monreal Gal\'{a}n, N., Pau, J.: Closure of Hardy spaces in the Bloch space. J. Math.
Anal. Appl. 429(2), 1214–-1221 (2015)
%
\bibitem{GKU1}
Gerg\"un, S.,   Kaptano\u glu, H. T.,   \"Ureyen, A. E.: Reproducing kernels for harmonic Besov spaces on the ball,
C. R. Math. Acad. Sci. Paris, 347, 735--738 (2009)
%
\bibitem{GKU2}
Gerg\"un, S.,   Kaptano\u glu, H. T.,   \"Ureyen, A. E.: Harmonic Besov spaces on the ball,
Int. J. Math., 27 (9), 1650070, 59 pp. (2016)
%
\bibitem{GZ}
Ghatage, P., Zheng, D.: Analytic functions of bounded mean oscillation and the Bloch space. Integral
Equ. Oper. Theory 17, 501–515 (1993)
%
\bibitem{GY}
G\"{o}\u{g}\"{u}\c{s}, N.G., Yılmaz, F.: Closures of Bergman–Besov spaces in the weighted Bloch spaces on the unit
ball. Complex Anal. Oper. Theory 15(6), 1-–13 (2021)
%
\bibitem{LR}
Liu, B., R\"{a}tty\"{a}, J.: Closure of Bergman and Dirichlet spaces in the Bloch norm. Ann. Acad. Sci. Fenn.
Math. 45(2), 771–-783 (2020)
%
\bibitem{MZ}
Manhas, J.S., Zhao, R.: Closures of Hardy and Hardy–Sobolev spaces in the Bloch type space on the
unit ball. Complex Anal. Oper. Theory 12(5), 1303–-1313 (2018)
%
\bibitem{M}
Miao, J.: Reproducing kernels for harmonic Bergman spaces of the unit ball,
Monatsh. Math., 125, 25--35 (1998)
%
\bibitem{MN}
Monreal Gal\'{a}n, N., Nicolau, A.: The closure of the Hardy space in the Bloch norm. Algebra i Analiz
22(1), 75–-81 (2010); reprinted in St. Petersburg Math. J. 22(1), 55–-59 (2011)
%
\bibitem{O1}
G.O. Okikiolu, Aspects of the theory of bounded  integral operators in $L^{p}$-spaces, Academic, London, 1971.
%
\bibitem{TM}
Tjani, M.: Distance of a Bloch function to the little Bloch space. Bull. Austral. Math. Soc. 74(1),
101-–119 (2006)
%
\bibitem{XW}
Xu, W.: Distances from Bloch functions to some Möbius invariant function spaces in the unit ball of
Cn. J. Funct. Spaces Appl. 7(1), 91-–104 (2009)
%
\bibitem{ZR}
Zhao, R.: Distances from Bloch functions to some Möbius invariant spaces. Ann. Acad. Sci. Fenn.
Math. 33(1), 303–-313 (2008)
\end{thebibliography}
\end{document}